\documentclass[12pt]{article}
\usepackage{latexsym}
\usepackage{amsfonts}
\def\C{{\bf C}}    

\def\R{{\bf R}}    
\def\N{{\bf N}}

\def\text{\hbox}   
\def\Im{\text{\rm Im~}}
\def\Re{\text{\rm Re~}}
\def\:{\colon}
\def\ds{\displaystyle\strut}
\newtheorem{thm}{Theorem}
\newtheorem{lem}[thm]{Lemma}
\newtheorem{cor}[thm]{Corollary}
\newtheorem{thmA}{Theorem A}

\newtheorem{thmB}{Theorem B}

\title{Zeros of differential polynomials in real
meromorphic functions}
\author
{Walter 
Bergweiler\thanks{Supported  
by the German-Israeli Foundation
for Scientific Research and Development (G.I.F.),
grant no.\ G-643-117.6/1999},
Alex Eremenko\thanks{Supported 
by NSF grants DMS-0100512 and DMS-0244421.}
\ and Jim Langley}
\date{June 30, 2004}
\begin{document}
\maketitle
\begin{abstract}
We investigate when differential polynomials
in real transcendental meromorphic functions have non-real zeros.
For example, we show  that if
$g$ is a real transcendental meromorphic function,
$c \in \R \setminus \{ 0 \}$
and $n\geq 3$ is an integer, then $g'g^n - c$ has infinitely many
non-real zeros. If $g$ has only finitely many poles, then
this holds for $n\geq 2$. Related results for rational functions
$g$ are also considered.
\end{abstract}

\section{Introduction and results}

Our starting point is the following result due to
Sheil-Small~\cite{ss} which solved a longstanding conjecture.
\begin{thmA}
Let $f$ be a real polynomial of degree $d$.
Then $f'+f^2$ has at least $d-1$ distinct non-real
zeros which are not zeros of $f$.
\end{thmA}

In the special case that
$f$ has only real roots this theorem is due to 
Pr\"ufer~\cite[Ch. V, 182]{ps};
see~\cite{ss} for further discussion of the result.

The following Theorem B is an analogue of Theorem A for
transcendental meromorphic functions. 
Here ``meromorphic'' will mean
``meromorphic in the complex plane'' unless explicitly stated otherwise.
A meromorphic function is called
real if it maps the real axis $\R$ to $\R\cup\{\infty\}$.
\begin{thmB}
Let $f$ be a real transcendental meromorphic
function with finitely many poles.
Then $f'+f^2$ has infinitely many non-real zeros 
which are not zeros of $f$.
\end{thmB}

Theorem B is a special case of~\cite[Theorem~1.3]{bel}.
The result that $f'+f^2$ has infinitely many non-real 
zeros if $f$ is a real entire transcendental function 
follows from the main theorem of~\cite{bf}.

A corollary of Theorem A is the following result.
\begin{cor} \label{thm1}
Let $f$ be a real polynomial of degree $d$
and $m\geq 2$ an integer.
Then $f'+f^m$ has at least $d(m-1)-1$ distinct non-real
zeros which are not zeros of $f$.
In particular, $f'+f^m$ has at most $d+1$ real zeros.
\end{cor}

In fact, putting $g(z)=f(w)^{m-1}$
where $w=z/(m-1)$ we have 
$$g'(z)+g(z)^2=f(w)^{m-2}\left(f'(w)+f(w)^m\right)$$
so that 
Corollary~\ref{thm1} follows from Theorem A applied to $g$,
since $g$ has degree $d(m-1)$.

The same argument yields 
the following corollary to Theorem B.
\begin{cor} \label{thm2}
Let $f$ be a real transcendental meromorphic
function with fi\-nite\-ly many poles 
and $m\geq 2$ an integer.
Then $f'+f^m$ has infinitely many non-real zeros 
which are not zeros of $f$.
\end{cor}

We note that Theorem B does not hold for meromorphic functions
with infinitely many poles,
a simple example being $f(z)=-\tan z$. Similarly Theorem A fails
for rational functions.

In this paper we consider to what extent Corollary~\ref{thm1} 
holds for rational functions, and Corollary~\ref{thm2} 
for meromorphic functions $f$ with infinitely many poles.
We recall that  Hayman \cite[Corollary to Theorem~9]{h} showed that if 
$f$ is a transcendental meromorphic function and $m\geq 5$,
then $f'+f^m$ has infinitely many zeros.
In fact, his proof shows that
$f'+f^m$ has infinitely many zeros which are not zeros of $f$.
(This can also be seen from the main result of~\cite{be} and 
formula (\ref{2}) below.)

Mues \cite{Mu} showed that Hayman's
result remains valid for $m=4$,
and the first and second author~\cite{be} proved that this also holds 
for $m=3$. We also note that if $f$ is rational and $f'+f^m$ has no
zeros for some $m\geq 3$, then $f$ is constant.

\begin{thm}\label{rat}
Let $f$ be a real rational function of degree $d$ and $m\geq 5$ an
integer.
Then $f'+f^m$ has at least $(m-4)d$ distinct non-real
zeros which are not zeros of $f$.
In particular,
 $f'+f^m$ has at most $4d$ real zeros.

If $m$ is odd, then 
$f'+f^m$ has at least $(m-3)d$ distinct non-real
zeros which are not zeros of $f$ so that
$f'+f^m$ has at most $3d$ real zeros.
\end{thm}

\begin{thm}\label{thm3}
Let $f$ be a real transcendental meromorphic function and $m\geq 5$
an integer.
Then $f'+f^m$ has infinitely many non-real zeros which are not zeros 
of~$f$.
\end{thm}
Theorem~\ref{thm3} improves 
Hayman's result in the case that $f$ is real. 
We will show by examples that the restriction $m\geq 5$ in this
theorem is best possible. 

If $g(z)=1/f(w)$ where $w=cz$, then
\begin{equation}\label{2}
f'(w)+f(w)^m=c^{-1}g(z)^{-m}(c-g(z)^{m-2}g'(z)),
\end{equation}
and we obtain the following 
results from Theorems~\ref{rat} and~\ref{thm3}.
\begin{cor}\label{crat}
Let $g$ be a real rational function of degree $d$ and $n\geq 3$
an integer.
Then for every real $c\neq 0$ the equation 
$g^ng'=c$ has at least $d(n-2)$ distinct non-real solutions.
\end{cor}
\begin{cor}\label{c1}
Let $g$ be a real transcendental meromorphic function and $n\geq 3$
an integer.
Then for every real $c\neq 0$ the equation $g^ng'=c$ has infinitely many
non-real solutions.
\end{cor}
The condition $n\geq 3$ is best possible 
in order to conclude that there are non-real solutions.
However for 
polynomials 
and transcendental meromorphic functions with finitely many poles we have
\begin{thm}\label{thm4pol}
Let $g$ be a real polynomial of degree $d$ and $n\geq 2$ an integer.
Then for every real $c\neq 0$ the number of
distinct non-real solutions of the equation 
$g^ng'=c$ is at least $d(n-1)-1$ if $n$ is even and at least
$d(n-1)$ if $n$ is odd.
In particular, this equation has at most $2d$ real solutions.
\end{thm}
\begin{thm}\label{thm4}
Let $g$ be a real transcendental meromorphic function with finitely
many poles and $n\geq 2$ an integer.
Then for every real $c\neq 0$ the equation $g^ng'=c$ has infinitely many
non-real solutions.
\end{thm}
Examples show that if $n=1$, then 
the equation $g^ng'=c$  need not have non-real solutions.

For our last result, we return to the value distribution of $f' + f^m$.
~Hayman~\cite{h} proved not only that if 
$f$ is a transcendental meromorphic
function and $m\geq 5$ then $f'+f^m$ has infinitely many zeros, 
but that under these conditions
$f'+f^m+c$ has infinitely many zeros for any $c\in\C$. Examples show
that when $m$ is $5$ or $6$ and $c \in \R \setminus \{0\}$ then
$f'+ f^m + c$ may fail to have non-real zeros, with $f$ real,
trancendental and meromorphic. 

\begin{thm}\label{cneq0}
Let $f$ be a real meromorphic function and set
$G = f' + f^m + c$, where 
$c \in \R \setminus \{ 0 \}$ and $m \in \N$.

Suppose that
$f$ is transcendental. If $m \geq 7$ then $G$ has infinitely many non-real
zeros which are not zeros of $f$. The same conclusion holds for $m \geq 4$
if $\overline{N}(r, f) = o( T(r, f))$ and thus in particular if $f$
has finitely many poles.

Suppose finally that
$f$ is a non-constant rational function with $f' \not \equiv - c$. If
$m \geq 6$ then $G$ has at least one non-real zero
which is not a zero of $f$, and the same conclusion holds for $m \geq 3$
if $f$ is a polynomial.
\end{thm}

\section{A result from complex dynamics}
One of our main tools is 
a result
from holomorphic dynamics.
Recall that the Fatou set of a non-linear
meromorphic function $f$ is
the set where the iterates $f^{\circ n}$ of $f$ are defined and form
a normal family.
We say that $\zeta\in\C$ is a multiple fixed point  of $f$
(of multiplicity~$\mu$)
if $f(z)-z$ has a multiple zero (of multiplicity~$\mu$) at $\zeta$.
This has to be slightly modified if $\zeta=\infty$.
We say that $\infty$ is a 
multiple fixed point  of $f$
(of multiplicity~$\mu$) 
if $1/f(1/z)$ has a multiple fixed point (of multiplicity~$\mu$) at $0$.
\begin{lem}\label{fatoutheorem}
Let $f$ be  a non-linear rational or 
transcendental meromorphic function
and let $\zeta$ be a multiple fixed point of $f$ of multiplicity
$\mu$.
Then there are $\nu:=\mu-1$ components 
$U_1,U_2,\dots,U_{\nu}$ of the Fatou set of
$f$ satisfying $f(U_j)\subset U_j$, $\zeta\in\partial U_j$ and
$f^{\circ n}(z)\to \zeta$ for all $z\in U_j$ as $n\to\infty$.
Each $U_j$ contains at least one singularity of 
$f^{-1}$, the inverse function of~$f$.

In addition, the $U_j$ can be labelled
such that if $z\in U_j$, then 
$${\rm arg}(f^{\circ n}(z)-\zeta)\to \theta_j: =\frac{-{\rm arg}\;f^{(\nu+1)}(\zeta)-\pi +2\pi j}{\nu}$$
as $n\to\infty$ if $\zeta\in\C$,
while if $\zeta=\infty$ then,
for some $\theta_0\in\R$,
$${\rm arg}(f^{\circ n}(z))\to \theta_j: =\frac{\theta_0 +2\pi j}{\nu}$$
as $n\to\infty$.
\end{lem}

The domains $U_j$ appearing in Lemma~\ref{fatoutheorem} 
are called {\em Leau domains}.
Note that the singularities of $f^{-1}$ are precisely the 
critical and asymptotic values of $f$.
If $f$ is rational, then  the only singularities of $f^{-1}$ 
are the critical values. In fact, the domains $U_j$ then
also contain critical points, but we do not need this result.

Lemma~\ref{fatoutheorem} is due to Fatou,
and can now be found in every textbook on complex dynamics.
Excellent introductions to complex dynamics are \cite{Mil,St};
see \cite[\S 10]{Mil} or
\cite[\S 3.5]{St} for a proof and discussion of the results stated 
in Lemma~\ref{fatoutheorem}.
We note that~\cite{Mil,St} and most other textbooks on complex 
dynamics treat only the case that $f$ is rational, but the 
proofs extend to the
case of transcendental meromorphic functions; see also~\cite[\S 4.3]{b}.

To demonstrate how dynamics
works we begin with a simple direct proof of Corollary~\ref{thm1}
found by the second author in 1989, after reading Sheil-Small's paper.
This dynamical proof is also reproduced in \cite{hink}.
\vspace{.1in}

{\em Proof of Corollary \ref{thm1}}.
Put 
\begin{equation}
\label{eq1}
F(z)=z-\frac{1}{f(z/(m-1))^{m-1}}.
\end{equation}
Then $F'(z)=0$ gives
$f'(w)+f(w)^m=0$
where $w=z/(m-1)$.
We have
$$F(z)=z+cz^{-d(m-1)}+O(z^{-d(m-1)-1}),\quad z\to\infty,$$
where $c\in\R\backslash\{0\}$.
This implies that $\infty$ is a multiple fixed point of 
$f$ of multiplicity $\mu=d(m-1)+2$.
Let $U_1,U_2,\dots,U_{\mu-1}$ be the Leau domains at
$\infty$ and let $\theta_1,\theta_2,\dots,\theta_{\mu-1}$
be as in Lemma~\ref{fatoutheorem}.
Each $U_j$ contains a critical value and, 
as $F$ is real, this critical value and the corresponding 
critical point can be real only if $\theta_j$ is a
multiple of $\pi$.
This is the case for at most $2$ of the $\mu-1$ values of $j$, and thus
$F$ has at least $\mu-3=d(m-1)-1$ non-real critical points.

As the total number of zeros of $f'+f^m$ is $dm$ we obtain the 
result.~\hfill$~\Box$
\vspace{.1in}

The 
proof of Theorem \ref{rat} will use the same idea.

\vspace{.1in}

{\em Proof of Theorem \ref{rat}}. 
Let $f(z)=P(z)/Q(z)$ where $P$ has degree $p$ and
$Q$ has degree $q$, so that $d=\max\{p,q\}$.
We consider again the function
$F$ defined by (\ref{eq1}).
If $q<p=d$, then 
$$F(z)=z+cz^{-(d-q)(m-1)}+O(z^{-(d-q)(m-1)-1}),\quad z\to\infty,$$
and the argument used in the proof of Corollary~\ref{thm1} shows that 
$F$ has a multiple fixed point at $\infty$ and that
the Leau domains associated to this fixed point
contain at least $(d-q)(m-1)-1$ non-real zeros of $F'$.

Let $\zeta_1,\dots,\zeta_k$ be the finite poles of $f$,
with multiplicities $s_1,\dots,s_k$.
Then $\zeta_j (m-1)$ is a fixed point of multiplicity $s_j(m-1)$
of $F$.
Lemma~\ref{fatoutheorem} now yields that 
of the $s_j(m-1)-1$ Leau domains associated to $\zeta_j$ at 
most $2$ can contain a real critical value, so 
at least $s_j(m-1)-3$ of these domains
give rise to non-real zeros of $F'$.

Overall this leads to $\sum_{j=1}^k s_j(m-1)-3=q(m-1)-3k$
non-real zeros of $F'$. 
If $q=d$ we thus have $d(m-1)-3k\geq d(m-1)-3d=d(m-4)$ non-real
zeros of $F'$.
If $q<d$, then we obtain
$(d-q)(m-1)-1+q(m-1)-3k=
d(m-1)-1-3k\geq d(m-1)-1-3(d-1)=d(m-4)+2$ non-real 
zeros of $F'$.

If $m$ is odd, then the number of Leau domains at 
the $\zeta_j$ and $\infty$ is odd. Lemma~\ref{fatoutheorem}
shows that each of them, with at most 1 exception,
contains a non-real critical value.
The above argument then shows that $F'$ has at least
$d(m-3)$ non-real zeros.~\hfill~$\Box$

\section{Proof of Theorem 4}
Our second main tool is the following result of Pang \cite{pang} which
already found many applications, see, for example, 
Zalcman's survey \cite{z}. 
\begin{lem}\label{l1}
Let $f$ be a meromorphic function with unbounded
spherical derivative, and 
$\kappa\in(-1,1)$. Then there exist sequences $z_j\in\C$ and $a_j>0$
such that 
\begin{equation}
\label{7}
a_j^{-\kappa}f(z_j+a_jz)\to h(z),\quad j\to\infty,
\end{equation}
uniformly on compact subsets of $\C$, where $h$ is a non-constant
meromorphic function with bounded spherical derivative.
Furthermore, one can choose
\begin{equation}
\label{tri}
z_j\to\infty\quad\mbox{and}\quad a_j\to 0,\quad j\to\infty.
\end{equation}
\end{lem}

The following result is an immediate consequence
of the definition of the order of a meromorphic function,
using the Ahlfors-Shimizu form of the Nevanlinna characteristic. 
\begin{lem}\label{l2} Let $f$ be a meromorphic function with 
bounded spherical derivative. Then $f$ is of order at most $2$.
\end{lem}

{\em Proof of Theorem~\ref{thm3}.}
The proof is by contradiction, assuming that
$f'+f^m$ has only finitely many non-real zeros which are not zeros of $f$.

The first step 
is to reduce the result to the case of functions with
bounded spherical derivative. 
Suppose $f$ has unbounded spherical derivative.
We apply Lemma~\ref{l1} with $\kappa=1/(1-m)$.
Then $\kappa m=\kappa-1$ and we obtain
$$a_j^{-\kappa m}f(z_j+a_jz)^m\to h(z)^m.$$
In $\C \setminus h^{-1} ( \{ \infty \} )$ 
we have
$$a_j^{-\kappa m}f'(z_j+a_jz)=a_j^{-\kappa+1}f'(z_j+a_jz)\to h'(z)$$
and
$$h'(z)+h(z)^m=\lim_{j\to\infty}a_j^{-\kappa m}
\left( f'(z_j+a_j z)+f(z_j+a_j z)^m\right).$$
If $|\Im z_j/a_j|\to\infty$ then all zeros of $h'+h^m$ are zeros of $h$. 
This contradicts the result
of Hayman \cite[Corollary to Theorem~9]{h} mentioned in 
the introduction; see also~\cite{be}. 

Thus 
$|\Im z_j/a_j|\not\to\infty$ and we may assume that
$\Im z_j/a_j\to s\in\R$. Putting 
$x_j=\Re z_j =z_j-i \;\Im z_j$ we find that 
$$a_j^{-\kappa}f(x_j+a_jz)=
a_j^{-\kappa}f(z_j+a_j(z-i\;\Im z_j/a_j))
\to h(z-is).$$
We may thus assume that $z_j\in\R$ so that $h$ is real,
since otherwise we can
replace $z_j$ by $x_j$ and $h(z)$ by $h(z-is)$.
We obtain a non-constant real meromorphic function $h$ with
bounded spherical derivative. It follows from (\ref{tri})
that {\em all} non-real zeros of $h'+h^m$ are zeros of $h$. 
Theorem~\ref{rat} implies that $h$ is transcendental.

If $f$ has bounded spherical derivative, then we put $h=f$.
Again $h$ is transcendental, but $h'+h^m$ may have finitely
many non-real zeros which are not zeros of $h$.

As in the proof of Corollary~\ref{thm1} and Theorem~\ref{rat}
we  consider the auxiliary function 
$$F(z)=z-\frac{1}{h(z/(m-1))^{m-1}},$$
and note that poles of $h$ give rise to
multiple fixed points of $F$ of multiplicity 
at least $m-1$.
Lemma~\ref{fatoutheorem} implies that at least $2$ of the 
(at least $3$) Leau 
domains associated to such a fixed point
contain a non-real singularity
of $F^{-1}$.

Theorem 1 from \cite{be} implies that all non-real asymptotic values
correspond to logarithmic singularities of the inverse function
$F^{-1}$. But according to the Denjoy--Carleman--Ahlfors Theorem
\cite[\S 258]{nev}, a function of order at most $2$ can have at
most $4$ logarithmic singularities. 
Thus $F$ has only finitely many non-real critical and
asymptotic values,
so we conclude that $h$ has finitely many poles.
We obtain a contradiction with Corollary~\ref{thm2}.~\hfill~$\Box$

\section{Proof of Theorems 7 and 8}
{\em Proof of Theorem~\ref{thm4pol}}.
Without loss of generality we can assume that $c=1$.
We consider the function
$$G(z)=z-\frac1{n+1} g(z)^{n+1}$$
and note that the critical points of $G$ are
exactly the solutions of $g(z)^ng'(z)=1$.

If $\zeta$ is a zero of $g$ 
of multiplicity $s$, then $\zeta$ is a multiple 
fixed point of $G$ of multiplicity $s(n+1)$. 
Lemma~\ref{fatoutheorem} yields that 
the $s(n+1)-1$ associated Leau domains 
contain at least $s(n+1)-3$ non-real critical values.
If $n$ is odd and thus the number
$s(n+1)-1$ of Leau domains is odd, then we even obtain
$s(n+1)-2$ non-real critical values in these domains.
Let $\zeta_1,\dots,\zeta_k$ be the
zeros of $g$, with multiplicities $s_1,\dots,s_k$.
If $n$ is odd then we
find that $G'$ has at least 
$(n+1)\sum_{j=1}^k s_j -2k=d(n+1)-2k$ non-real zeros.
Since $k\leq d$ this implies that $G'$ has at least
$d(n-1)$ non-real zeros.

Now we consider the case that $n$ is even.
If $\zeta$ is a simple zero of $g$, then the previous argument
based on Lemma~\ref{fatoutheorem} will only yield that $n-2$ 
of the $n$ associated Leau domains contain a non-real
critical value, which does not give any information 
in the most interesting case that $n=2$.  We note, however, 
that if $g'(\zeta)<0$, then $G^{(n+1)}(\zeta)=-n! g'(\zeta)^{n+1}>0$.
With $U_j$ as in 
Lemma~\ref{fatoutheorem} we find that if $z\in U_j$, then
$${\rm arg}(G^{\circ \ell}(z)-\zeta)\to 
\theta_j: =\frac{-{\rm arg}\;G^{(n+1)}(\zeta) -\pi+2\pi j}{n}
=\frac{(2j-1)\pi}{n}$$
as $\ell\to\infty$.
As $n$ is even, $(2j-1)/n$ is never an integer and thus
all $n$ Leau domains at $\zeta$ contain a non-real critical value
in this case.

Let $K$ be the number of real simple zeros $\zeta$ 
of $g$ for which $g'(\zeta)>0$.
Between two such zeros there must be a real simple zero $\zeta$ 
satisfying $g'(\zeta)<0$ or a multiple zero of odd multiplicity.
Let $L$ be the number of real simple zeros $\zeta$ 
of $g$ for which $g'(\zeta)<0$ and let 
$M$ be the number of real multiple zeros of odd multiplicity.
Then 
\begin{equation}
\label{L}
L+M+1\geq K.
\end{equation}
We denote by $N$ the total number of multiple zeros of $g$
and by $P$ the number of non-real simple zeros.
Thus $K+L+P$ is the total number of simple zeros of $g$.
We find that 
\begin{equation}
\label{d}
d\geq K+L+P+3M+2(N-M)=K+L+P+M+2N.
\end{equation}
Denote by $\xi_1,\dots,\xi_N$ the multiple zeros of $g$, 
with multiplicities $t_1,\dots,t_N$.
The above considerations based on Lemma~\ref{fatoutheorem}
show that there are at least $(n+1)t_j-3$ non-real critical
values of $G$ contained in the Leau domains associated to
$\xi_j$. 
Also, for a non-real simple zero of $g$ each of the 
$n$ associated Leau domains contains a non-real
critical value of $G$.
Overall the number $\nu$ of non-real critical
points of $G$ thus satisfies  
\begin{eqnarray*}
\nu 
&\geq &
(n+1)\sum_{j=1}^N t_j -3N +K(n-2)+Ln+Pn\\
&=&
(n+1)(d-K-L-P) -3N+K(n-2)+Ln+Pn\\
&=&
(n-1)d +2d-3N-3K-L-P.
\end{eqnarray*}
Using (\ref{d}) and (\ref{L}) we obtain
\begin{eqnarray*}
\nu 
&\geq &
(n-1)d +2(K+L+P+M+2N)-3N-3K-L-P\\
&=&
(n-1)d -K+L+2M+N+P\\
&\geq &
(n-1)d -(L+M+1)+L+2M+N+P\\
&=&
(n-1)d -1+M+N+P.
\end{eqnarray*}
The conclusion follows since $M,N,P\geq 0$.~\hfill~$\Box$
\vspace{.1in}

To prove Theorem~\ref{thm4}
we cannot use the reduction to functions of finite order based on
Lemma~7, because there exists a non-constant real entire function,
namely $h(z)=z$, with the property that all solutions of the equation
$h'(z)h(z)^2=1$ are real. We need instead a direct proof that $g$
is of finite order. 
\vspace{.1in}

{\em Proof of Theorem~\ref{thm4}}. In view of Corollary~\ref{c1}
it is enough to consider the case $n=2$. We can assume without loss
of generality that $c=1$.

Suppose that the equation $g'(z)g(z)^2=1$ has finitely many non-real solutions.
First we show that $g$ has order
at most $1$. To accomplish this, we need
the characteristic function in the upper half-plane as developed
by Tsuji \cite{Tsuji0} and Levin and Ostrovskii \cite{LeO} (see also \cite{GO}),
and as used in \cite{bel}. 
For $\psi$ meromorphic and non-constant in the closed upper half-plane
${\rm Im} \, z \geq 0$ and for  
$t \geq 1$ let
$\mathfrak{n} (t, \psi) $ be the number of poles of $\psi$, counting
multiplicity, in $\{ z:|z-it/2|\leq t/2, |z| \geq 1\}$, and set
$$
\mathfrak{N} (r, \psi) = \int_1^r
\frac{ \mathfrak{n} (t, \psi) }{t^2} dt, \quad r \geq 1.
$$
The Tsuji characteristic is 
$$
\mathfrak{T} (r, \psi)  =
\mathfrak{m} (r, \psi)  +
\mathfrak{N} (r, \psi),
$$
where
$$\mathfrak{m} (r, \psi) =
\frac1{2 \pi} \int_{ \sin^{-1} (1/r)}^{ \pi - \sin^{-1} (1/r)}
\frac{ \log^+ | \psi(r \sin \theta e^{i \theta } )|}
{r \sin^2 \theta } d \theta.$$
We refer the reader to \cite{bel,GO,LeO,Tsuji0} for the fundamental
properties of the Tsuji characteristic, but note in particular that
the lemma on the logarithmic derivative
\cite[p.332]{LeO} (see also \cite[Theorem 3.2, p.141]{GO}) gives
\begin{equation}
\mathfrak{m} (r, \psi'/\psi) = O(\log r + \log^+ \mathfrak{T} (r, \psi) )
\label{t2}
\end{equation}
as $r \to \infty $ outside a set of finite measure.

The following lemma is a direct
analogue for $g$ of a result of 
Hayman from \cite{h}.

\begin{lem}\label{lemt2}
We have
$$
\mathfrak{T} (r, g) = O( \log r ), \quad r \to \infty .
$$
\end{lem}

{\em Proof.}
We follow Hayman's proof as in \cite{h}, but using the Tsuji characteristic
and in particular (\ref{t2}).
Let 
$$
\phi (z) = \frac13 g(z)^3 .
$$
Then $\phi $ has finitely many poles, and $\phi ' - 1$ has
finitely many zeros in the open upper half-plane $H$, and so
\begin{equation}
\mathfrak{N} \left( r, \phi  \right) 
+
\mathfrak{N} \left( r, \frac1{\phi ' - 1} \right)  = O(1) .
\label{milloux1}
\end{equation}
Milloux' inequality
\cite[Theorem 3.2, p.57]{Haym2} translates directly in terms
of the Tsuji characteristic
to give, using (\ref{t2}) and (\ref{milloux1}), outside a set
of finite measure,
\begin{equation}
\mathfrak{T}(r, \phi) < 
\mathfrak{N} \left( r, \frac1{\phi} \right) 
-
 \mathfrak{N}_0 \left( r, \frac1{\phi''} \right) 
+ O(\log r + \log^+\mathfrak{T} (r, \phi))
\label{ts1}
\end{equation}
in which $\ds\mathfrak{N}_0 \left( r, \frac1{\phi''} \right)$ counts only
zeros of $\phi ''$ which are not multiple zeros of $\phi' - 1$.
But all zeros of $\phi $ have multiplicity at least 3 and so are zeros
of $\phi ''$ but not zeros of $\phi' - 1$, and consequently each such
zero contributes 2 to
$\ds\mathfrak{n} \left( r, \frac1{\phi} \right) -
\mathfrak{n}_0 \left( r, \frac1{\phi''} \right) $ but at least 3 to
$\ds\mathfrak{n} \left( r, \frac1{\phi} \right)$.
Hence (\ref{ts1}) becomes
\begin{eqnarray*}
\mathfrak{T}(r, \phi)& <& 
\frac23 \mathfrak{N} \left( r, \frac1{\phi} \right)
+ O(\log r + \log^+\mathfrak{T} (r, \phi))\\
&<&
\frac23 \mathfrak{T}(r, \phi)
+ O(\log r + \log^+\mathfrak{T} (r, \phi)).
\end{eqnarray*} 
Thus
$3 \mathfrak{T} (r, g  ) \leq 
\mathfrak{T} (r, \phi  ) + O( \log r) = O( \log r )$ initially outside
a set of finite measure, and hence without exceptional set since
$\mathfrak{T} (r, g  ) $ differs from a non-decreasing function
by a bounded term \cite{Tsuji0} (see 
also \cite[p.980]{bel}).~\hfill~$\Box$

\begin{lem}\label{lolem}
The function $g$ has order at most $1$.
\end{lem}

{\em Proof.}
This
proof is almost identical to \cite[Lemma 3.2]{bel}
and to corresponding arguments in \cite{LeO}.
Lemma \ref{lemt2}
 and an inequality of Levin-Ostrovskii 
\cite[p.332]{LeO} (see also \cite[Lemma 2.2]{bel}) give
\begin{equation}
\int_R^\infty \frac{m_{0\pi}(r, g)}{r^3} dr
\leq  \int_R^\infty \frac{\mathfrak{m} (r, g)}{r^2} dr
=O\left(\frac{\log R}{R}\right),
\quad R  \to \infty,
\label{ts2}
\end{equation}
in which
$$
m_{0\pi} (r, g) = 
\frac1{2 \pi} \int_0^\pi \log^+ | g(r e^{i \theta} ) | d \theta .
$$
But $g$ 
is real on the real axis and has finitely many poles and so 
$$
T(r, g) = m(r, g)  + O( \log r )  = 2 m_{0\pi} (r, g ) 
+ O( \log r )
$$
and (\ref{ts2}) now gives, as $R  \to \infty$,
$$
\frac{ T(R, g)}{ R^{2}}
\leq 2 \int_R^\infty \frac{T(r, g)}{r^3} dr
\leq  4 \int_R^\infty \frac{\mathfrak{m} (r, g)}{r^2} dr 
+O\left(\frac{\log R}{R^2}\right)
=O\left(\frac{\log R}{R}\right),
$$
from which the lemma follows.~\hfill~$\Box$
\vspace{0.1in}

The rest of the proof is similar to the proof of Theorem~\ref{thm4pol}.
However, the arguments simplify considerably since
we have restricted attention to the case $n=2$ and since we do not have
to count the number of non-real critical points as precisely as in the
proof of Theorem~\ref{thm4pol}.

We put 
$$G(z)=z-\frac13 g(z)^3.$$
Then the zeros of $G'$ are exactly the solutions
of $g(z)^2g'(z)=1$,
so all critical points of $G$, with finitely many exceptions, are real.
By Lemma~\ref{lolem} and the Denjoy-Carleman-Ahlfors theorem,
$G$ has at most $2$ finite asymptotic values.
Thus the number of non-real singularities of $G^{-1}$ is finite.

The fixed points of $G$ are all multiple
and they coincide with zeros of $g$.
If $g$ 
has finitely many zeros then $g(z)=p(z)\exp(az),$ where
$p$ is a rational function and $a\in\R\backslash\{0\}$. For such $g$,
it is easy to see that the equation
$g'(z)g(z)^2=1$ has infinitely many non-real solutions.

Hence we may assume that $g$ has infinitely many zeros.
A non-real zero of $g$ is a non-real multiple fixed point
of $G$, and Lemma~\ref{fatoutheorem}
implies that the Leau domains associated to it 
contain a non-real singularity of $G^{-1}$.
A multiple zero of $g$ is a multiple fixed point of $G$ of 
multiplicity at least $6$, and Lemma~\ref{fatoutheorem}
yields that at least $3$ of the associated Leau domains 
contain a non-real singularity of $G^{-1}$. As the number
of non-real singularities of $G^{-1}$ is finite, we see
that only finitely many zeros of $g$ are non-real or multiple,
and thus all but finitely many of 
the zeros of $g$ are real and simple.
This implies that 
there are infinitely many real zeros $\zeta$ of $g$ with
$g'(\zeta)<0$.
As in the proof of Theorem~\ref{thm4pol} we see that
there are at least $2$ non-real singularities of $G^{-1}$ 
contained in the Leau domains associated to such a point $\zeta$.
Since the number
of non-real singularities of $G^{-1}$ is finite, we 
deduce that there are only finitely many 
real zeros $\zeta$  of $g$ 
satisfying
$g'(\zeta)<0$, 
a contradiction.~\hfill~$\Box$
\vspace{.1in}

\section{Proof of Theorem \ref{cneq0}}

Let $m \geq 3$, let $f$, $G$ and $c$ be as in the
hypotheses, and define $g$ and $H$ by
$$
 g = \frac{f' + c}{f^m }, \quad
\frac{f'+c}{G} = \frac{g}{g+1}  , 
$$
\begin{equation}
f'' - \left( \frac{G'}{G} \right)  (f' + c) =
\frac{g'G}{(g+1)^2} = \frac{ f^m g'}{g+1} = f^m H .
\label{tc0}
\end{equation}
The third equation in (\ref{tc0}) is obtained
by differentiating the second and multiplying
by $G$.
The function $g$ is non-constant, 
poles of $g$ are zeros of $f$, and $g(z) = -1$ implies $G(z) = 0$.
We may assume that $H \not \equiv 0$, since 
$H \equiv 0$ implies that
$g$ is constant.
Let $S(r, f)$ denote any quantity which 
is $o( T(r, f) )$ as $r \to \infty$ outside a set of finite measure.

\begin{lem}
We have, as $r \to \infty $,
\begin{equation}
(m-1) m(r, f) = m(r, f^{m-1}) \leq m \left(r, \frac{1}{H} \right)  
- \alpha \log r + S(r, f),
\label{tc11}
\end{equation}
in which: $\alpha = 0$ if $f$ is transcendental;
$\alpha = 2$ if $f( \infty ) = \infty$;
$\alpha = 2$ if $f(\infty) \in \C \setminus \{ 0 \}$ and $g( \infty ) \neq -1$;
$\alpha = 1$ otherwise.
\end{lem}

{\em Proof.} For transcendental $f$ we apply the method of Clunie's lemma 
\cite[Lemma 3.3, p.68]{Haym2}. Separating into the two cases
$|f| < 1, |f| \geq 1$ and using (\ref{tc0}) we obtain
$$
|f^{m-1}H| \leq |H| + 
\left| \frac{f''}{f} \right| + 
\left| \frac{G'}{G} \right|\left( \left| \frac{f'}{f} \right| + |c| \right) , 
\quad m(r, f^{m-1} H) = S(r, f),
$$
and (\ref{tc11}) follows on writing $f^{m-1} =
(f^{m-1}H)/H$. Suppose now that $f$ is rational.
If $f$ has
a pole of multiplicity $p \in \N$ at $\infty$, 
then $g$ has a zero of multiplicity $p(m-1) + 1$ and
$H$ a zero of multiplicity $p(m-1) + 2$ at $\infty$, which gives
(\ref{tc11}) with $\alpha = 2$. 
If $f( \infty ) \in \C$ then $m(r, f) = O(1)$ and
$H$ has at least a simple zero at $\infty$, 
while if $f(\infty) \in
\C \setminus \{ 0 \}$ and $g( \infty ) \neq -1$ then 
the zero of $H$ at $\infty$ is at least double.~\hfill~$\Box$
\vspace{.1in}

From Jensen's formula and (\ref{tc11}) we obtain at once, since
$m(r, H) = S(r, f)$,
\begin{equation}
(m-1) T(r, f) \leq (m-1) N(r, f) - N \left(r, \frac{1}{H} \right)  + N(r, H) 
- \alpha \log r + S(r, f).
\label{tc1}
\end{equation}
Now if $f$ has a pole in $\C$ of multiplicity $q$ then $g$ has a zero of 
multiplicity $mq - q - 1 \geq 1$ and so
$H$ has a zero of
multiplicity $mq - q - 2 = (m-1)q - 2$. So the contribution
to $(m-1)n(r, f) - n(r, 1/H)$ from this pole is $ 2$. 
A pole of $H$ is simple, and must be a zero or pole of $g + 1$, and
poles of $g$ are zeros of $f$. Thus
(\ref{tc1}) yields
\begin{eqnarray}
(m-1) T(r, f) &\leq&  2 \overline{N}(r, f) 
- N_0 \left( r, \frac{g+1}{g'} \right)  
+ \overline{N} \left( r, \frac1{f} \right) + \nonumber\\
&& + \overline{N}_1 \left( r, \frac1{g+1} \right)  
- \alpha \log r + S(r, f).
\label{tc2}
\end{eqnarray}
Here
$\ds N_0 \left( r, \frac{g+1}{g'} \right) $
counts zeros of $\ds \frac{g'}{g+1}$ which are not poles of $f$, and 
$\ds \overline{N}_1 \left( r, \frac1{g+1} \right)$ counts zeros of $g+1$
which are not zeros of $f$.
Write
$$
\overline{N}_1 \left( r, \frac1{g+1} \right)
= \overline{N}_{1,R} \left( r, \frac1{g+1} \right)
+ \overline{N}_{1,NR} \left( r, \frac1{g+1} \right),
$$
in which the subscripts $R, NR$ denote real, non-real zeros respectively.

\begin{lem}
We have, as $r \to \infty$,
\begin{eqnarray}
\overline{N}_{1,R} \left( r, \frac1{g+1} \right)
&\leq & \overline{N}(r, f) + 
\overline{N} \left( r, \frac1{f} \right) + 
N_0 \left( r, \frac{g+1}{g'} \right) + \nonumber\\
&& + \beta \log r  + S(r, f),
\label{tc3}
\end{eqnarray}
in which: $\beta = 0$ if $f$ is transcendental;
$\beta = -1$ if $f$ is rational and $g(\infty) = -1$; 
$\beta = 1$ if $f$ is rational and $g(\infty) \neq -1$.
\end{lem}

{\em Proof.}
Applying Rolle's theorem shows that between adjacent real zeros 
$x_1, x_2$ of $g+1$
there must be at least one point $x_0$
which is a pole of $g$ or a zero of $g'$ but
not a zero of $g+1$. In the first case $x_0$
is a zero of $f$, while in the the second case $x_0$ either
is a pole of $f$ or contributes to
$\ds N_0 \left( r, \frac{g+1}{g'} \right) $. If $g(\infty) = -1$ then the
same argument may be applied on intervals of form
$(-\infty, x_1), (x_2, \infty)$.
~\hfill~$\Box$
\vspace{.1in}

Combining (\ref{tc2}) and (\ref{tc3})
yields
\begin{eqnarray}
(m-1) T(r, f) &\leq&  3 \overline{N}(r, f) + 2
\overline{N} \left( r, \frac1{f} \right) +
\overline{N}_{1,NR} \left( r, \frac1{g+1} \right) + \nonumber\\ 
&& + (\beta - \alpha  ) \log r + S(r, f).
\label{tc4}
\end{eqnarray}
For transcendental $f$ both assertions of the theorem follow at once from
(\ref{tc4}).

Suppose now that $f$ is a rational function,
and that all non-real zeros of $g+1$ are zeros of $f$.
If $f(\infty) = \infty$ then $\beta - \alpha \leq -1$ and a 
contradiction arises from (\ref{tc4}) if
$m \geq 6$ or if $m \geq 3$ and $f$ is a polynomial.
Next, if $f(\infty) = 0$ then $\beta - \alpha \leq 0$ and so
(\ref{tc4}) gives a contradiction for $m \geq 6$ since 
$\ds \overline{N} \left( r, \frac1{f} \right) < (1 - \varepsilon )
T(r, f)$ as $r \to \infty$, for some $\varepsilon > 0$.
Finally, if $f(\infty) \in \C \setminus \{ 0 \}$ then
$\beta - \alpha \leq -1$ and again 
if $m \geq 6$ we obtain a contradiction from (\ref{tc4}).
~\hfill~$\Box$
\vspace{.1in}

\section{Examples}

\noindent
1. $f(z)=e^{-z}+1$ gives $f'+f=1$ so one cannot take $m=1$ in
Theorem~\ref{thm3}.
\vspace{.1in}

\noindent
2. $f(z)=-\tan z$ gives $f'+f^2=-1$, so one cannot take $m=2$ in
Theorem~\ref{thm3}.
\vspace{.1in}

\noindent
3. 
$f(z)=1/(2\sin z)$ gives
$$f'(z)+f(z)^3=\frac{1-2\sin(2z)}{8\sin^3 z}$$
and this has only real zeros, so Theorem~\ref{thm3} does 
not hold with $m=3$.

$g(z)=1/f(z)=2\sin z$ gives $g'(z)g(z)-1=2\sin2z-1$ which has only real zeros.
So one cannot put $n=1$ in Theorem~\ref{thm4}.

\vspace{.1in}

\noindent
4. For $f(z)=a-\tan a^4z$ with $a=12$ we obtain
$$w(t)=f'+f^4=-a^4(1+t^2)+(a+t)^4,\quad\mbox{where}\quad t=-\tan a^4z  .$$
Indeed, $w(0)=0,\; w'(0)>0$ and $w(\pm\infty)=+\infty$.
It follows that $w$ has a negative root.
Now $w(1)=-2\times 12^4+13^4=-12911<0$ which implies that $w$ has $2$
positive roots. Thus all $4$ roots of $w$ are real. The full
preimage of the real line under $-\tan a^4z$ is contained in
the real line, so all roots of $f'+f^4$ are real.
So Theorem~\ref{thm3} does not hold with $m=4$.

Another example with this property is given by
\[
f(z)= \frac {\sqrt{2} + \tan(4x)}{1 + \sqrt{2}\tan(4x)}
\]
so that
\[
f'(z)+f(z)^4=-\frac{\tan^3 (4x) \left(4\sqrt{2}+
7\tan (4x)\right)}{\left(1 + \sqrt{2}\tan(4x)\right)^4}.
\]

\vspace{.1in}

Taking $g=1/f$ we see that Corollary~\ref{c1} does not hold with
$n=2$.
\vspace{.1in}

\noindent
5. For $f(z)=1/z$ the function
$f'(z)+f(z)^4=(1 - z^2)/z^4$ has only real zeros,
so the the condition $m\geq 5$ in Theorem \ref{rat}
cannot be weakened to $m\geq 4$ in order to conclude that 
$f'+f^m$ has non-real zeros.
Equivalently $g(z)=1/f(z)=z$
shows that in Corollary \ref{crat}  
the condition $n\geq 2$ does not yield the existence of
non-real solutions of $g'g^n=1$.

An example of degree $2$ with these properties is given by
$$f(z)=\frac{5z+3}{2z^2+8z+3}$$
and $g(z)=1/f(z)$.
Here $0$ is a zero of $g'g^2-1$ of multiplicity $3$, and
the other $3$ zeros of this function are also real.
This can be seen by numerical computation, but can also
be deduced from Lemma~\ref{fatoutheorem}, since  $G(z)=z-g(z)^3/3$
has $4$ Leau domains 
associated to the $2$
real zeros of $g$, each of them containing a critical point, and symmetric
about the real axis since $g' \geq 0$ on $\R$.
If one contained a non-real 
critical point, then 
it would also contain the complex conjugate of this
point, and this leads to a contradiction since there are only
$4$ distinct critical points.

\vspace{.1in}

\noindent
6. The following examples shows that for transcendental meromorphic $f$
one cannot take $m= 5$ or $m = 6$ in Theorem~\ref{cneq0}. Let 
$f(z)=-\tan z$. Then
$$f'(z)+f(z)^6+1= (\tan z+1)(\tan z-1)(\tan^2 z +1)\tan^2 z$$
has only real zeros. 

For the case $m=5$ let
\begin{equation}
5a^4 - 10a^2 + 1 = 0, \quad
b = 5a - 10 a^3 , \quad c = - a^5 - b .
\label{aeqn}
\end{equation}
Set $f(z) = a + \tan (bz) = a + t$. Then
$$
G = f' + f^5 + c = t (t^2 + 1)( t^2 + 5 a t + 10 a^2 - 1).
$$
Equation (\ref{aeqn}) has a root $a =\sqrt{1-2/\sqrt{5}}
= 0 . 3249 \ldots $ for which we have
$25 a^2 - 4 ( 10 a^2 - 1 )  > 0$, so that all zeros of $G$ are real.

The example
$$
f(z) = - 16 z^2 + 8z + 2, \quad
f'(z) + f(z)^2 - 12 = 256 z^3 (z - 1),
$$
shows that Theorem~\ref{cneq0} is sharp for polynomial $f$.

We suspect that
Theorem~\ref{cneq0} is not sharp for non-polynomial rational functions,
nor for transcendental functions with finitely many poles.


{\em W.B.:
Mathematisches Seminar,

Christian--Albrechts--Universit\"at zu Kiel,

Ludewig--Meyn--Str.\ 4,

D--24098 Kiel,

Germany

{bergweiler@math.uni-kiel.de}

\vspace{.1in}

A.E.: Purdue University,

West Lafayette IN 47907

USA

eremenko@math.purdue.edu
\vspace{.1in}

J.K.L.: School of Mathematical Sciences,

University of Nottingham,

Nottingham NG7 2RD

UK

jkl@maths.nottingham.ac.uk}
\end{document}